\documentclass{amsart}
 \makeatletter
 \def\section{\@startsection{section}{1}%
 	\z@{.7\linespacing\@plus\linespacing}{.5\linespacing}%
 	{\bfseries
 		\centering
 	}}
 	\def\@secnumfont{\bfseries}
 	\makeatother
\newtheorem{thm}{Theorem}[section]
\newtheorem{lem}{Lemma}[section]

\newtheorem{prp}{Proposition}[section]
\theoremstyle{remark}
\newtheorem{rem}{Remark}[section]
\theoremstyle{definition}

\newcommand{\real}{\mathbb{R}}
\newcommand{\sig}{\mathcal{F}}

\begin{document}

\title
[On the Solution of SFDEs via Memory Gap]
{On the Solution of Stochastic Functional Differential Equations via Memory Gap}

\author{Flavia Sancier}
\address{Flavia Sancier: Sciences Division,  Antioch College, 1 Morgan Place, Yellow Springs, OH 45387, USA.}
\email{fsancier@antiochcollege.edu}

\author{Salah Mohammed}
\address{Salah Mohamed: Deceased December 21, 2016}



\begin{abstract}
We present an alternative proof for the existence of solutions of stochastic functional differential equations satisfying a global Lipschitz condition. The proof is based on an approximation scheme in which the continuous  path dependence does not go up to the present: there is a \emph{memory gap}. Strong convergence is obtained by \emph{closing the gap}. Such approximation is particularly useful when extending stochastic models with discrete delay to models with continuous full finite memory.
\end{abstract}

\maketitle

\section{Introduction}

Stochastic systems with memory are processes whose evolution in time is governed by random forces as well as an intrinsic dependence of the state on a finite part of its past history. Such systems are described by stochastic functional differential equations (SFDEs) \cite{mohammed}. The existence and uniqueness of a general class of SFDEs was first derived by Mohammed in \cite{mohammed} (Chapter II, sec. 2) through a successive approximation scheme, and conditions on the coefficients have been improved by others (e.g. \cite{von}).

Our results are based upon approximating stochastic systems with full finite memory by stochastic systems with a \emph{memory gap}. 
Solutions of systems with a memory gap are processes in which the intrinsic dependence of the state on its history goes only up to a specific time in the past. In this way, there is a \emph{gap} between the past and present states.

The paper is outlined as follows. In section \ref{frmwk}, we introduce the framework necessary for the development of our results. In section \ref{apsch}, we describe an approximation scheme and show the properties necessary to obtain convergence. In section \ref{cgap}, we show strong convergence of the approximation scheme to the solution of an SFDE with full finite memory.

\section{Framework}\label{frmwk}
Let $|\cdot|$ stand for the ($d$-dimensional) Euclidean norm, and let $|\cdot|_{d\times m}$ stand for the Frobenius norm in the space of real $d\times m$ matrices, $\real^{d\times m}$. 
Denote by $C_d:=C([-L,0],\real^d)$ the Banach space of all continuous paths $\eta:[-L,0]\rightarrow \real^d$ given the supremum norm $$\|\eta\|_{C_d}:=\sup_{v\in[-L,0]}|\eta(v)|.$$ 

Consider a filtered probability space  $(\Omega, \mathcal{F}, (\mathcal{F}_t)_{t\in[0,T]}, P)$ satisfying the usual conditions, viz. the filtration $(\sig_t)_{t\in[0,T]}$ is right-continuous, and each $\sig_t$, $t\in[0,T]$, contains all $P$-null sets in $\sig$. For a Banach space $E$, denote by $L^2(\Omega,E)$ the Banach space of all (equivalent classes of) $(\sig,$Borel$E)$-measurable maps $\Omega\rightarrow E$ which are $L^2$ (in the Bochner sense). The norm in $L^2(\Omega,E)$ is given by
$$\|\eta\|_{L^2(\Omega,E)}:=\left[\int_{\Omega}\|\eta(\omega)\|_E^2dP(\omega)\right]^{1/2},$$
for any $\eta\in {L^2(\Omega,E)}$.

Let $F:[0,T]\times L^2(\Omega,C_d)\rightarrow L^2(\Omega,\real^d)$ and $G:[0,T]\times L^2(\Omega,C_d)\rightarrow L^2(\Omega,\real^{d\times m})$ be jointly continuous and uniformly Lipschitz in the second variable, viz.
\begin{eqnarray}\label{lipschitz}
&&\nonumber\|F(t,\psi_1)-F(t,\psi_2)\|_{L^2(\Omega,\real^d)}+\|G(t,\psi_1)-G(t,\psi_2)\|_{L^2(\Omega,\real^{d\times m})}\\
&&\leq \alpha \|\psi_1-\psi_2\|_{L^2(\Omega,C_d)}
\end{eqnarray}
for all $t\in[0,T]$ and $\psi_1, \psi_2 \in L^2(\Omega,C_d)$. The Lipschitz constant $\alpha$ is independent of $t\in[0,T]$. 
For each $(\mathcal{F}_t)_{t\in[0,T]}$-adapted process $y:[0,T]\rightarrow L^2(\Omega,C_d)$, the processes $F(\cdot,y(\cdot))$ and $G(\cdot,y(\cdot))$ are also $(\mathcal{F}_t)_{t\in[0,T]}$-adapted. 
Let $W$ be an $m$-dimensional Brownian Motion on the filtered probability space $(\Omega, \mathcal{F}, (\mathcal{F}_t)_{t\in[0,T]}, P)$ and let $\theta\in L^2(\Omega, C_d)$ be an $\mathcal{F}_0$-measurable process.
\begin{rem}
	The functionals $F$ and $G$ also satisfy the linear growth property 
	\begin{equation}\label{lingrowth}
	\|F(t,\psi)\|_{L^2(\Omega,\real^d)}+\|G(t,\psi)\|_{L^2(\Omega,\real^{d\times m})}\leq D(1+\|\psi\|_{L^2(\Omega,C_d)}),
	\end{equation}
	where $D$ is a positive constant independent of $t\in [0,T]$. To see this, set $\psi_1=\psi$ and $\psi_2=0$ in the Lipschitz condition (\ref{lipschitz}), and use the joint continuity of $F$ and $G$.
\end{rem}
Consider the stochastic functional differential equation
\begin{equation}\label{sfde}
\left\{ \begin{array}{ll}
dx(t)=F(t,x_t)dt+G(t,x_t)dW(t), \quad t\in [0,T]\\
x(t)= \theta(t), \quad t\in [-L,0],
\end{array}\right.
\end{equation}
where the segment $x_t\in L^2(\Omega,C_d)$ is given by the relation $x_t(\omega)(s):=x(\omega)(t+s)$, $s\in[-L,0]$, a.a. $\omega\in\Omega$.
A solution of (\ref{sfde}) is a process $x\in L^2(\Omega,C([-L,T],\real^d))$ adapted to $(\mathcal{F}_t)_{t\in[0,T]}$, with initial process $\theta$, which satisfies the It\^{o} integral equation
\begin{equation*}
x(t)=\left\{ \begin{array}{ll}
\theta(\cdot)(0)+\int_0^t{F(u,x_{u})du}+\int_0^t{G(u,x_{u})dW(u)}, \quad t\in [0,T]\\
\theta(\cdot)(t), \quad t\in [-L,0],
\end{array}\right.
\end{equation*}
almost surely. In order to prove existence of solutions of the SFDE (\ref{sfde}), we will make frequent use of 
a martingale-type inequality for the Ito integral, 
which we state below. 

%
\begin{lem}\label{mohisometry}
	Let $W:[a,b]\times\Omega\rightarrow\real^m$ be an $m$-dimensional Brownian Motion on a filtered probability space $(\Omega, \mathcal{F}, (\mathcal{F}_t)_{t\in[a,b]}, P)$. Suppose $g:[a,b]\times\Omega\rightarrow L(\real^m,\real^d)$ is measurable, $(\mathcal{F}_t)_{t\in[a,b]}$-adapted and $\int_a^b E|g(t,\cdot)|^{2k}dt<\infty$, for a positive integer $k\geq 1$. Then
	$$ E\sup_{t\in[a,b]}\left|\int_a^t g(u,\cdot)dW(u)\right|^{2k}\leq A_k(b-a)^{k-1}\int_a^b E|g(u,\cdot)|^{2k}du,$$
	where $$A_k:=d^{k-1}\left(\frac{4k^3m^2}{2k-1}\right)^k.$$
\end{lem}
\noindent \emph{Proof.} For a proof, the reader may refer to Mohammed \cite{mohammed} (pg. 27).



\section{The approximation scheme}\label{apsch}
In this section, we introduce an approximation scheme for the linear SFDE (\ref{sfde}), and show the properties necessary to obtain convergence in $L^2(\Omega, C([-1-L,T],\real^d))$.  

Let $\hat{\theta}$ be given by
\begin{equation}\label{hattheta}
\hat{\theta}(t):=\left\{ \begin{array}{ll}
\theta(t), \quad t\in [-L,0],\\
\theta(-L), \quad t\in [-1-L,-L],
\end{array}\right.
\end{equation}
and define the sequence of processes, for any positive integer $k\geq 1,$
\begin{equation}\label{sfdek}
x^k(t)=\left\{ \begin{array}{ll}
\theta(0)+\int_0^t{F(u,x_{u-1/k}^k)du}+\int_0^t{G(u,x_{u-1/k}^k)dW(u)}, \quad t\in [0,T]\\
\hat{\theta}(t), \quad t\in [-1-L,0],
\end{array}\right.
\end{equation}
where for $t\in[-1,0]$, $x_t^k$ is given by $x^k_t(s)=\hat{\theta}_t(s):=\hat{\theta}(t+s)$, $s\in [-L,0]$, a.a. $\omega\in \Omega$, $k\geq 1$. For $t\in[-1,0]$ define $\sig_t:=\sig_0$. In what follows, we show existence of the $x^k$'s by integrating forward over steps of length $1/k$.



\begin{prp}\label{step1}
	Each $x^k$ satisfies the following properties for any $t\in[-1/k,T]$:\\
	\noindent (i) $x^k(t)$ and $x^k_t$ are well defined and $\mathcal{F}_t$-measurable.\\
	\noindent (ii) $\left(x^k(u)\right)_{u\in[-1/k,t]}\in L^2(\Omega, C([-1/k,t],\real^d))$ and $x_{t}^k\in L^2(\Omega,C_d)$, with
	$$E \left[\sup_{v\in[-1/k,t]}|x^k(v)|^2\right] + \|x_{t}^k\|_{L^2(\Omega,C_d)}^2\leq K,$$ where $K$ is a constant independent of $k$.
\end{prp}
\noindent\emph{Proof.} We use induction with forward steps of length $1/k$. For simplicity, consider $T$ a positive integer. For $t\in[-1/k,0]$, properties $(i)$ and $(ii)$ hold trivially. 
Now assume that properties (i) and (ii) hold for any $t\in[(n-1)/k,n/k]$, where $n$ is an integer satisfying $0\leq n<kT$. 
Then we show that properties (i) and (ii) hold for any $t\in [n/k,(n+1)/k]$. 
For any $t\in [n/k,(n+1)/k]$, one can write 
\begin{equation}\label{induck}
x^k(t)=x^k(n/k)+\int_{n/k}^t{F(u,x_{u-1/k}^k)du}+\int_{n/k}^t{G(u,x_{u-1/k}^k)dW(u)}.
\end{equation}
Since $\left(x^k(t)\right)_{t\in[(n-1)/k,n/k]}$ is well defined, $(\sig_t)_{t\in[(n-1)/k,n/k]}$-adapted and continuous (induction hypothesis), then $(x^k_t)_{t\in[(n-1)/k,n/k]}$ is well defined,\\ $(\sig_t)_{t\in[(n-1)/k,n/k]}$-adapted and continuous (Lemma II-2.1 in Mohammed \cite{mohammed}). This implies that $x^k_{t-1/k}$ is $\sig_{t-1/k}$-measurable for any $t\in [n/k,(n+1)/k]$. Therefore, by the assumptions in section \ref{frmwk} and the previous discussion, the processes $[n/k,(n+1)/k]\ni t\mapsto F(t,x^k_{t-1/k})$ and $[n/k,(n+1)/k]\ni t\mapsto G(t,x^k_{t-1/k})$ are continuous and $(\sig_t)_{t\in[n/k,(n+1)/k]}$-adapted. Hence, from (\ref{induck}), $\left(x^k(t)\right)_{t\in[n/k,(n+1)/k]}$ is a well defined $(\sig_t)_{t\in[n/k,(n+1)/k]}$-adapted continuous semimartingale. This shows that property (i) holds for any $t\in[n/k,(n+1)/k]$. 

Next, notice that for any $t\in [-1/k,(n+1)/k]$,
\begin{eqnarray}\label{indkn}
\nonumber E \left[ \sup_{s\in [-L,0]}|x^k_{t}(s)|^2 \right]\leq  \|\theta\|^2_{L^2(\Omega,C_d)}+E \left[ \sup_{s\in [0,t]}|x^k(s)|^2\right].
\end{eqnarray}
Now applying lemma \ref{mohisometry} and the linear growth property (\ref{lingrowth}) of $F$ and $G$, we have for any $t\in [0,(n+1)/k]:$
{\allowdisplaybreaks
	\begin{eqnarray}\label{calc1}
	& &E\left [ \sup_{v\in[0,t]}|x^k(\cdot)(v)|^2\right ]\nonumber \\
	&=&E\left [  \sup_{v\in[0,t]}\left| \theta(0)+\int_0^v{F(u,x^k_{u-1/k})du}+\int_0^v{G(u,x^k_{u-1/k})dW(u)} \right|^2 \right ]\nonumber \\
	&\leq&E\left [ 3|\theta(0)|^2\right]+3t\left [ \int_0^{t}{\|F(u,x^k_{u-1/k})\|_{L^2(\Omega,\real^d)}^2du}\right]\nonumber \\
	&+&3(4m^2)\left[\int_0^{t}\|{G(u,x^k_{u-1/k})\|_{L^2(\Omega,\real^{d\times m})}^2du}\right]\nonumber\\
	&\leq& 3\|\theta\|_{L^2(\Omega,C_d)}^2+6TD^2(T+4m^2)(1+\|\theta\|_{L^2(\Omega,C_d)}^2)\nonumber\\
	&+&6D^2(T+4m^2)\int_0^{t}E\left [ \sup_{v\in[0,u]}|x^k(v)|^2\right ]du.
	\end{eqnarray}
}
Applying Gronwall's inequality to the above inequality, we obtain:
\begin{eqnarray*}
	E \left[\sup_{v\in[0,t]}|x^k(v)|^2\right]\leq C_1e^{C_2T}
\end{eqnarray*}
where $C_1:=3\|\theta\|_{L^2(\Omega,C_d)}^2+6TD^2(T+4m^2)(1+\|\theta\|_{L^2(\Omega,C_d)}^2)$ and $C_2:=6D^2(T+4m^2)$. Hence, for each $t\in[0,(n+1)/k]$,
\begin{eqnarray*}
	E \left[\sup_{v\in[-1/k,t]}|x^k(v)|^2\right]+\|x^k_{t}\|_{L^2(\Omega,C)}\leq 2C_1e^{C_2T} +2\|\theta\|_{L^2(\Omega,C)}=:K.
\end{eqnarray*}
Notice that the constant $K$ does not depend on $k$ or $t$.
This concludes the induction argument.$\quad_\square$


\begin{prp}\label{step2}
	For any integer $\gamma\geq 1$, each $x^k$ satisfies $E|x^k(t)-x^k(s)|^{2\gamma}\leq B_\gamma|t-s|^\gamma$ for all $s,t\in[0,T]$, where $B_\gamma$ is a constant independent of $k$.
\end{prp}
\noindent \emph{Proof.} By lemma \ref{mohisometry}, proposition \ref{step1} (iii), and the linear growth property of $F$ and $G$ (\ref{lingrowth}), we obtain for any $0\leq s<t\leq T$:
{\allowdisplaybreaks
	\begin{eqnarray*}
		& & \hspace*{-6mm}E\left[\left|x^k(t)-x^k(s)\right|^{2\gamma}\right]=E\left[\left|\int_s^tF(u,x^k_{u-1/k})du+\int_s^tG(u,x^k_{u-1/k})dW(u)\right|^{2\gamma}\right]\\
		&\leq&2^{2\gamma-1}|t-s|^{2\gamma-1}\int_s^t\|F(u,x^k_{u-1/k})\|^{2\gamma}_{L^2(\Omega,\real^d)}du\\
		&+&2^{2\gamma-1}A_\gamma|t-s|^{\gamma-1}\int_s^t\|G(u,x^k_{u-1/k})\|^{2\gamma}_{L^2(\Omega,\real^{d\times m})}du\\
		&\leq& 2^{2\gamma-1}D^{2\gamma} (T^\gamma+A_\gamma)|t-s|^{\gamma-1}(1+\sqrt{K})^{2\gamma}|t-s|=B_\gamma|t-s|^{\gamma},
	\end{eqnarray*}
}
where $A_\gamma:=d^{\gamma-1}(\frac{4\gamma^3 m^2}{2\gamma-1})^\gamma$ and $B_\gamma:=2^{2\gamma-1}D^{2\gamma} (T^\gamma+A_\gamma)(1+\sqrt{K})^{2\gamma}$.
By symmetry in $t$ and $s$, the proposition holds for any $t,s\in [0,T]$. $\quad_\square$

Next, we state Kolmogorov's continuity criterion for a sequence of Banach-valued stochastic processes. The theorem will be used in the proof of proposition \ref{step3}.

\begin{thm}\label{kolm}
	\textbf{Kolmogorov's continuity criterion for a sequence of stochastic processes}. Let $\{X^k(t)\}_{k=1}^\infty$, $t\in[0,T]$, be a sequence of stochastic processes with values in a Banach space $E$. Assume that there exist positive constants $\rho_1$, $c$ and $\rho_2>1$, all independent of $k$, satisfying
	\begin{equation}\label{unifbound}
	E[\|X^k(t)-X^k(s)\|_E^{\rho_1}]\leq c|t-s|^{\rho_2},
	\end{equation}
	for every $s,t\in [0,T]$. Then each $X^k$ has a continuous modification $\tilde{X}^k$.
	Further, let $b$ be an arbitrary positive number less than $\frac{\rho_2-1}{\rho_1}$. Then there exists a positive random variable $\xi_k$ with $E[\xi_k^{\rho_1}]<H$, where H is a constant independent of $k$, such that
	\[
	\|\tilde{X}^k(t)-\tilde{X}^k(s)\|_E\leq \xi_k|t-s|^b,
	\]
	for every $s,t\in [0,T]$ and a.s..
\end{thm}
\noindent \emph{Proof.} The reader may refer to Kunita \cite{kunita}, pg. 31, for a proof. Although the author does not consider a sequence of stochastic processes satisfying (\ref{unifbound}), it is easy to check that the theorem indeed holds for such a sequence.


\begin{prp}\label{step3}
	Let $\beta\in(0,1/2)$ be a fixed constant. Each $x^k$ satisfies
	\begin{eqnarray*}
		&(i)&\,\,|x^k(t)-x^k(s)|\leq c_k|t-s|^\beta\,\,\mbox{ for all } s,t\in[0,T]\,\, a.s.;\\
		&(ii)&\,\,\|x^k_t-x^k_s\|^2_{L^2(\Omega,C_d)}\leq 3\tilde{c}|t-s|^{2\beta}+\\
		&& 2E\sup_{v\in(-(t\wedge L)\wedge0,-(s\wedge L)\wedge0]}|\hat{\theta}(0)-\hat{\theta}(s+v)|^2+\\
		&& E\sup_{v\in[-L,-(t\wedge L)\wedge0]}|\hat{\theta}(t+v)-\hat{\theta}(s+v)|^2\,\,\mbox{ for all } -1\leq s<t\leq T,\,\, a.s.,
	\end{eqnarray*}
	where $\tilde{c}$ is a constant independent of $k$ and $c_k$ is a positive random variable satisfying $E(c_k^\gamma)\leq \tilde{c}$.
\end{prp}
\noindent \emph{Proof.}
Let $\rho>\frac{1}{1-2\beta}$ be an integer. From proposition \ref{step2}, $E|x^k(t)-x^k(s)|^{2\rho}\leq B_\rho|t-s|^{\rho}$, for any $s,t\in[0,T]$. Since $\beta<\frac{\rho-1}{2\rho}$, then it follows from Kolmogorov's continuity criterion (theorem \ref{kolm}) that there exists a positive random variable $c_k$ such that $|x^k(t)-x^k(s)|\leq c_k|t-s|^\beta$ a.s., with $E(c_k^\gamma)\leq \tilde{c}$, where $\tilde{c}$ is a constant independent of $k$. This proves part (i).

We now proceed to prove part (ii). For any $-1\leq s<t\leq T$, 
\begin{eqnarray}
\nonumber \|x^k_t-x^k_s\|_{L^2(\Omega,C_d)}^2 &=& E\sup_{v\in[-L,0]}|x^k(t+v)-x^k(s+v)|^2\\ \label{es1}
&\leq&E\sup_{v\in(-(s\wedge L)\wedge0,0]}|x^k(t+v)-x^k(s+v)|^2\\ \label{es2}
&+&E\sup_{v\in(-(t\wedge L)\wedge0,-(s\wedge L)\wedge0]}|x^k(t+v)-x^k(s+v)|^2\hspace*{5mm} \\ \label{es3}
&+&E\sup_{v\in[-L,-(t\wedge L)\wedge0]}|x^k(t+v)-x^k(s+v)|^2.
\end{eqnarray}
Using part (i), we obtain the estimate for  (\ref{es1}):
\begin{eqnarray}
\nonumber& &E\sup_{v\in(-(s\wedge L)\wedge0,0]}|x^k(t+v)-x^k(s+v)|^2\leq E \sup_{v\in(-(s\wedge L)\wedge0,0]}\left(c_k|t-s|^\beta\right)^2\\
\nonumber &=&E(c_k^2|t-s|^{2\beta})=E(c_k^2)|t-s|^{2\beta}\leq \tilde{c}|t-s|^{2\beta}.
\end{eqnarray}
Moreover, from part (i), we also estimate (\ref{es2}):
\begin{eqnarray}
\nonumber& &E\sup_{v\in(-(t\wedge L)\wedge0,-(s\wedge L)\wedge0]}|x^k(t+v)-x^k(s+v)|^2\\
\nonumber&\leq& 2E\sup_{v\in(-(t\wedge L)\wedge0,-(s\wedge L)\wedge0]}\{|x^k(t+v)-x^k(0)|^2+|\hat{\theta}(0)-\hat{\theta}(s+v)|^2\}\\
\nonumber&\leq& 2E\sup_{v\in(-(t\wedge L)\wedge0,-(s\wedge L)\wedge0]}\{c_k^2|t+v|^{2\beta}+|\hat{\theta}(0)-\hat{\theta}(s+v)|^2\}\\
\nonumber&\leq& 2\tilde{c}|t-s|^{2\beta}+2E\sup_{v\in(-(t\wedge L)\wedge0,-(s\wedge L)\wedge0]}|\hat{\theta}(0)-\hat{\theta}(s+v)|^2.
\end{eqnarray}
Finally, (\ref{es3}) becomes:
\begin{eqnarray}
\nonumber E\sup_{v\in[-L,-(t\wedge L)\wedge0]}|x^k(t+v)-x^k(s+v)|^2=E\sup_{v\in[-L,-(t\wedge L)\wedge0]}|\hat{\theta}(t+v)-\hat{\theta}(s+v)|^2.
\end{eqnarray}
Hence, for $-1\leq s<t\leq T,$
\begin{eqnarray*}
	\nonumber \|x^k_t-x^k_s\|_{L^2(\Omega,C_d)}^2 &\leq&3\tilde{c}|t-s|^{2\beta}+2E\sup_{v\in(-(t\wedge L)\wedge0,-(s\wedge L)\wedge0]}|\hat{\theta}(0)-\hat{\theta}(s+v)|^2\\
	&+&E\sup_{v\in[-L,-(t\wedge L)\wedge0]}|\hat{\theta}(t+v)-\hat{\theta}(s+v)|^2. \quad_\square
\end{eqnarray*}
Proposition \ref{step3}(i) shows that each $x^k$ is pathwise $\beta$-H\"{o}lder continuous on [0,T]. Notice that if $t>s>L$, part (ii) can be written as $\|x^k_t-x^k_s\|^2_{L^2(\Omega,C_d)}\leq 3\tilde{c}|t-s|^{2\beta}$, which implies that the memory process $t\mapsto x^k_{t}$ is $\beta$-H\"{o}lder continuous (as an $L^2(\Omega,C)$-valued function) only for $t\in[L,T]$. On the other hand, if the initial process $\theta$ is pathwise $\beta$-H\"{o}lder continuous, then $t\mapsto x^k_{t}$ is $\beta$-H\"{o}lder continuous (as an $L^2(\Omega,C)$-valued function) for all $t\in[0,T]$. We  prove this in theorem \ref{rate}.


\section{Closing the gap}\label{cgap}
In this section, we show that the sequence $(x^k)_{k=1}^\infty$ converges to the solution of (\ref{sfde}). The proof is divided into the following steps:\\
\noindent 1. The sequence $(x^k)_{k=1}^\infty$ converges to a limit $x\in L^2(\Omega,C([-1-L,T],\real^d))$. \\
\noindent 2. The process $\left(x(t)\right)_{t\in[-L,T]}$ satisfies the SFDE (\ref{sfde}).\\
\noindent 3. The process $\left(x(t)\right)_{t\in[-L,T]}$ is the unique solution of the SFDE (\ref{sfde}).\\
\noindent 4. If the initial process $\theta$ is pathwise $\beta$-H\"{o}lder continuous for a fixed $\beta\in(0,1/2)$, then the rate at which $x^k$ converges to $x$ is $\beta$.
\begin{prp}\label{step4}
	The sequence $(x^k)_{k=1}^\infty$ converges to a limit $x\in L^2(\Omega,C([-1-L,T],\real^d))$. 
\end{prp}
\noindent\emph{Proof}. We first notice that for any $t\in [-1,T]$,
\begin{eqnarray}\label{trajlk}
\nonumber\|x^l_t-x^k_t\|_{L^2(\Omega,C_d)}^2
&\leq&E \sup_{v\in[-1-L,t]}|x^l(v)-x^k(v)|^2=E \sup_{v\in[0,t]}|x^l(v)-x^k(v)|^2.\quad\quad
\end{eqnarray}
Now, from lemma \ref{mohisometry}, the Lipschitz condition (\ref{lipschitz}), and inequality (\ref{trajlk}), we have that for any $t\in[0,T]$ and $l>k$:
{\allowdisplaybreaks
	\begin{eqnarray}\label{convstep4}
	& &E \sup_{v\in[0,t]}|x^l(v)-x^k(v)|^2=E \sup_{v\in[0,t]}\left|\int_0^v[F(u,x^l_{u-1/l})-F(u,x^k_{u-1/k})]du \right.\\
	\nonumber& &+\left.\int_0^v[G(u,x^l_{u-1/l})-G(u,x^k_{u-1/k})]dW(u) \right|^2\\
	\nonumber&\leq&2t\int_0^t\|F(u,x^l_{u-1/l})-F(u,x^k_{u-1/k})\|_{L^2(\Omega,\real^d)}^2du\\
	\nonumber&& +2(4m^2)\int_0^t\|G(u,x^l_{u-1/l})-G(u,x^k_{u-1/k})\|_{L^2(\Omega,\real^{d\times m})}^2du\\
	\nonumber&\leq&4\alpha^2(t+4m^2)\int_0^t\left(\|x^l_{u-1/l}-x^l_{u-1/k}\|^2_{L^2(\Omega,C_d)}+\|x^l_{u-1/k}-x^k_{u-1/k}\|^2_{L^2(\Omega,C_d)}\right)du\\
	\nonumber&\leq&4\alpha^2(T+4m^2)\int_0^t\|x^l_{u-1/l}-x^l_{u-1/k}\|^2_{L^2(\Omega,C_d)}du\\
	\nonumber& &+4\alpha^2(T+4m^2)\int_0^tE \sup_{v\in[0,u]}|x^l(v)-x^k(v)|^2du.
	\end{eqnarray}
}
Hence, by Gronwall's inequality,
\begin{eqnarray}\label{gronwall1}
&&\nonumber E \sup_{v\in[0,t]}|x^l(v)-x^k(v)|^2\leq\\ &&4\alpha^2(T+4m^2)e^{4\alpha^2(T+4m^2)t}\int_0^t\|x^l_{u-1/l}-x^l_{u-1/k}\|^2_{L^2(\Omega,C)}du. 
\end{eqnarray}
From proposition \ref{step3}(ii) and the uniform continuity of $\hat{\theta}$, it follows that 
$\int_0^t\|x^l_{u-1/l}-x^l_{u-1/k}\|^2_{L^2(\Omega,C_d)}du\rightarrow 0$ as $l,k\rightarrow\infty$. Therefore, from inequality (\ref{gronwall1}), 
$$E \sup_{v\in[0,t]}|x^l(v)-x^k(v)|^2\rightarrow 0\quad\mbox{as}\quad l,k\rightarrow\infty.$$
This shows that the sequence $(x^k)_{k=1}^\infty$ is a Cauchy sequence in $L^2(\Omega,C([-1-L,T],\real^d))$ and therefore convergent to a limit $x\in L^2(\Omega,C([-1-L,T],\real^d))$. From (\ref{trajlk}), it also follows that for each $t\in[-1,T]$, $(x^k_t)_{k=1}^{\infty}$ converges to $x_t$ in $L^2(\Omega,C_d).\quad_\square$


\begin{prp}\label{step5}
	The process $\left(x(t)\right)_{t\in[-L,T]}$ satisfies the SFDE (\ref{sfde}).
\end{prp}
\noindent \emph{Proof.} To show this, we take limits as $k\rightarrow\infty$ in both sides of (\ref{sfdek}). The left-hand side of (\ref{sfdek}) converges to $x$ in $L^2(\Omega,C([-1-L,T],\real^d))$. Furthermore, $\left(x(t)\right)_{t\in[0,T]}$ is $(\sig_t)_{t\in[0,T]}$-adapted, since each $\left(x^k(t)\right)_{t\in[0,T]}$ is. Moreover, in a calculation similar to (\ref{convstep4}), we have that
\begin{eqnarray}\label{inconv}
\nonumber&E& \sup_{v\in[0,t]}\left|\int_0^v[F(u,x_{u})-F(u,x^k_{u-1/k})]du + \int_0^v[G(u,x_{u})-G(u,x^k_{u-1/k})]dW(u) \right|^2\\
\nonumber&\leq&4\alpha^2(t+4m^2)\int_0^t\|x_{u}-x_{u-1/k}\|^2_{L^2(\Omega,C_d)}du\\
& &+4\alpha^2(t+4m^2)\|x_{u-1/k}-x^k_{u-1/k}\|^2_{L^2(\Omega,C_d)}du.
\end{eqnarray}
The continuity of $[-1,T]\ni t\mapsto x_t\in L^2(\Omega,C_d)$ implies that
\begin{eqnarray*}
	\int_0^t\|x_{u}-x_{u-1/k}\|^2_{L^2(\Omega,C_d)}du\rightarrow 0 \,\, \mathrm{as} \,\,k\rightarrow\infty.
\end{eqnarray*}
Also, 
\begin{eqnarray*}
	\|x_{u-1/k}-x^k_{u-1/k}\|_{L^2(\Omega,C_d)}^2
	\leq E\sup_{v\in[-L,T]}|x(v)-x^k(v)|^2
	\rightarrow 0 \,\, \mathrm{as} \,\,k\rightarrow\infty.
\end{eqnarray*}

Hence, (\ref{inconv}) converges to 0 as $k\rightarrow\infty$.
This shows that the right-hand side of (\ref{sfdek}) converges to $\theta(0)+\int_0^{(\cdot)}{F(u,x_{u})du}+\int_0^{(\cdot)}{G(u,x_{u})dW(u)}$ in $L^2(\Omega,C([-L,T],\real^d))$ as $k\rightarrow\infty$. Therefore, $\left(x(t)\right)_{t\in[-L,T]}$ satisfies the sfde (\ref{sfde}).$\quad_\square$


\begin{prp}\label{step6}
	(Uniqueness) If $\tilde{x}$ is an $(\sig_t)_{t\in[0,T]}$-adapted process satisfying (\ref{sfde}), then $\tilde{x}=x|_{[-L,T]}$ a.s..
\end{prp}
\noindent \emph{Proof.} In a calculation similar to (\ref{convstep4}), we find the difference
\begin{eqnarray*}\label{uniq1}
	\nonumber & &\|\tilde{x}-x|_{[-L,T]}\|_{L^2(\Omega,C([-L,T],\real))}^2 = E\sup_{v\in[0,T]}|\tilde{x}(v)-x(v)|^2\\
	\nonumber&\leq&4\alpha^2(T+4m^2)\int_0^t\|\tilde{x}_{u}-x_{u}\|^2_{L^2(\Omega,C_d)}du\\
	& &+4\alpha^2(T+4m^2)\int_0^tE \sup_{v\in[0,u]}|\tilde{x}(v)-x(v)|^2du\\
	&\leq&8\alpha^2(T+4m^2)\int_0^tE \sup_{v\in[0,u]}|\tilde{x}(v)-x(v)|^2du.
\end{eqnarray*}
Hence, from Gronwall's inequality, it follows that 
\[
E\sup_{v\in[0,T]}|\tilde{x}(v)-x(v)|^2=0\Rightarrow \tilde{x}=x|_{[-L,T]}\,\,\, a.s..\quad_\square
\]

Thus far, we have proven existence and uniqueness of a strong solution $x$ to the SFDE (\ref{sfde}). Moreover, we have proven that the SFDE (\ref{sfde}) can be approximated by the sequence $(x^k)_{k=1}^\infty$ described in (\ref{sfdek}). The following theorem gives the rate in which this sequence converges to $x$, when the initial process satisfies a H\"{o}lder continuity  condition.
\begin{thm}\label{rate}
	Let $\beta\in(0,1/2)$ be a fixed constant. If the initial process $\theta$ satisfies 
	\begin{equation}\label{ahol}
	E|\theta(t)-\theta(s)|^{2\gamma}\leq C_{\theta}|t-s|^\gamma,
	\end{equation}
	for any $\gamma>1$, where $C_\theta$ is a positive constant, then $\theta$ is pathwise $\beta$-H\"{o}lder continuous and 
	\begin{equation*}
	E \sup_{v\in[0,t]}|x(v)-x^k(v)|^2\leq c\left(\frac{1}{k}\right)^{2\beta},
	\end{equation*}
	where c is a constant independent of $k$.
\end{thm}
\noindent \emph{Proof.} 
Let $\rho>\frac{1}{1-2\beta}$ be an integer. From (\ref{ahol}) and the fact that $\beta<\frac{\rho-1}{2\rho}$, Kolmogorov's continuity criterion (theorem \ref{kolm}) implies that there exists a positive random variable $c_\theta$ such that $|\theta(t)-\theta(s)|\leq c_\theta|t-s|^\beta$ a.s., with $E(c_\theta^\gamma)\leq \tilde{c_\theta}$, where $\tilde{c_\theta}$ is a positive constant. That is, $\theta$ is pathwise $\beta$-H\"{o}lder continuous.

Then notice that $\hat{\theta}$ is also pathwise $\beta$-H\"{o}lder continuous. Indeed, for $-1-L\leq s<t\leq 0$, 
\begin{eqnarray*}
	|\hat{\theta}(t)-\hat{\theta}(s)|=\left\{\begin{array}{lll}
		|\theta(t)-\theta(s)|\leq c_\theta|t-s|^{\beta}, \quad -L\leq s,t \leq 0\\
		|\theta(t)-\theta(-L)|\leq c_\theta|t+L|^{\beta}\leq c_\theta|t-s|^{\beta}, \quad s<-L,\,\,t>-L;\\
		|\theta(-L)-\theta(-L)|=0\leq c_\theta|t-s|^{\beta},\quad -1-L\leq s,t <-L.
	\end{array}
	\right.
\end{eqnarray*}
Then proposition \ref{step3}(ii) and the pathwise $\beta$-H\"{o}lder continuity of $\hat{\theta}$ imply that 
\begin{eqnarray}
\nonumber \|x^k_t-x^k_s\|^2_{L^2(\Omega,C_d)}&\leq& 3\tilde{c}|t-s|^{2\beta}+2E\sup_{v\in(-(t\wedge L)\wedge0,-(s\wedge L)\wedge0]}c_\theta^2|s+v|^{2\beta}\\
\nonumber&+&E\sup_{v\in[-L,-(t\wedge L)\wedge0]}c_\theta^2|t-s|^{2\beta}\\
\nonumber&\leq&3\tilde{c}|t-s|^{2\beta}+2E(c_\theta^2)|t-s|^{2\beta}+E(c_\theta^2)|t-s|^{2\beta}\\
\nonumber&\leq& 3(\tilde{c}+\tilde{c_\theta})|t-s|^{2\beta},
\end{eqnarray}
for any $s,t\in[-1,T]$.
Hence, for $l>k>0$,
\begin{eqnarray}
\nonumber\|x^l_{u-1/l}-x^l_{u-1/k}\|^2_{L^2(\Omega,C_d)}&\leq&3(\tilde{c}+\tilde{c_\theta})|1/k-1/l|^{2\beta},
\end{eqnarray}
which implies that  $\int_0^t\|x^l_{u-1/l}-x^l_{u-1/k}\|^2_{L^2(\Omega,C_d)}du\leq3T(\tilde{c}+\tilde{c_\theta})|1/k-1/l|^{2\beta}$, $t\in[0,T]$.
Therefore, from inequality \ref{gronwall1}, we obtain
\begin{equation*}
E \sup_{v\in[0,t]}|x^l(v)-x^k(v)|^2\leq 4\alpha^2(T+4m^2)e^{4\alpha^2(T+4m^2)t}3(\tilde{c}+\tilde{c_\theta})T|1/k-1/l|^{2\beta}.
\end{equation*}
Finally, letting $l\rightarrow\infty$ and $c:=12T\alpha^2(\tilde{c}+\tilde{c_\theta})(T+4m^2)e^{4\alpha^2(T+4m^2)T}$, we obtain
\begin{equation*}
E \sup_{v\in[0,t]}|x(v)-x^k(v)|^2\leq c\left(\frac{1}{k}\right)^{2\beta},\quad t\in[0,T].\quad_\square
\end{equation*}


\section{Remarks}\label{remarks}

\subsection{Order of convergence}
Theorem \ref{rate} gives the order of convergence for the approximation scheme (\ref{sfdek}), when the initial process $\theta$ is pathwise $\beta$-H\"{o}lder continuous. Since the quantity $1/k$ is the length of a small time interval, the order of convergence is given with respect to time increments. This result is comparable with the $0.5$ order of convergence of the Strong Euler-Maruyama scheme for SFDEs investigated in \cite{Hu,Ahmed}. 

\subsection{Existence of SODEs}
If we set $L=0$ in (\ref{sfde}) and (\ref{sfdek}), the approximation scheme reduces to the one in \cite{Bell}, which provides an alternative existence theorem for (non-delayed) stochastic ordinary differential equations (SODEs).

\subsection{Alternative scheme}
Alternatively, one could use the approximation scheme
\begin{equation}\label{sfdek2}
x^k(t)=\theta(0)+\int_0^t{F(u-1/k,x_{u-1/k}^k)du}+\int_0^t{G(u-1/k,x_{u-1/k}^k)dW(u)},
\end{equation}
for $t\in[0,T]$, and $x^k(t)=\hat{\theta}(t),$ for $t\in [-1-L,0]$. The process $\hat{\theta}$ is given by (\ref{hattheta}). In this case, the functionals $F$ and $G$ need to satisfy the additional regularity condition:
\begin{eqnarray}\label{reg}
\nonumber& &\|F(t,\eta)-F(s,\eta)\|_{L^2(\Omega,\real^d)}+\|G(t,\eta)-G(s,\eta)\|_{L^2(\Omega,\real^{d\times m})}\\
&\leq& \alpha_1(1+\|\eta\|_{C_d})|t-s|,
\end{eqnarray}
for any $\eta\in C_d$, $s,t\in[0,T]$, where $\alpha_1$ is a positive constant, independent of $\eta$, $s$ and $t$.

\subsection{Numerical simulation}
The approximation scheme (\ref{sfdek}) provides a way of numerically simulating the SFDE (\ref{sfde}), with the help a numerical method that can approximate Wiener integrals. More specifically, for each $n/k<t\leq (n+1)/k$ with $0\leq n<kT$, $x^k$ can be written as
\begin{equation*}
x^k(t)=x^k(n/k)+\int_{n/k}^t{F(u,x_{u-1/k}^k)du}+\int_{n/k}^t{G(u,x_{u-1/k}^k)dW(u)}. 
\end{equation*}
The integral $\int_{n/k}^t{F(u,x_{u-1/k}^k)du}$ is a Riemann integral and can be easily approximated. The integral $\int_{n/k}^t{G(u,x_{u-1/k}^k)dW(u)}$ is a Wiener integral, since $G(u,x_{u-1/k}^k)$ is $\sig_{n/k}$-measurable for each $u\in(n/k,(n+1)/k]$. Wiener integrals can be approximated in the Riemann-Stieltjes sense. In \cite{kloeden}, the authors provide several numerical schemes for SODEs, which include the less general case of a Wiener integral.


\end{document}